\documentclass{article}
\usepackage{amssymb}
\usepackage{amsfonts}
\usepackage{amsmath}

\newtheorem{theorem}{Theorem}

\newtheorem{lemma}[theorem]{Lemma}

\begin{document}

\title{Dunkl generalization of Sz\'{a}sz Beta type operators}
\author{Bayram \c{C}ekim$^{a}$,\"{U}lk\"{u} Dinlemez$^{b}$, \.{I}smet Y\"{u}%
ksel$^{c}$ \and Gazi University, Faculty of Science, Department of
Mathematics \and Teknikokullar Ankara/ TURKEY \and $^{a}$e-mail:
bayramcekim@gazi.edu.tr \and $^{b}$e-mail: ulku@gazi.edu.tr \and $^{c}$%
e-mail: iyuksel@gazi.edu.tr}
\maketitle

\begin{abstract}
The goal in the paper is to advertise Dunkl extension of Sz\'{a}sz beta type
operators. We initiate approximation features via acknowledged Korovkin and
weighted Korovkin theorem and obtain the convergence rate from the point of
modulus of continuity, second order modulus of continuity, the Lipschitz
class functions, Peetre's $K$-functional and modulus of weighted continuity
by Dunkl generalization of Sz\'{a}sz beta type operators.


\end{abstract}

\section{Introduction}

Newly, several mathematicians have made many studies concerning
generalization of Sz\'{a}sz operators (for example, see \cite{Atakut,
Ata-Buyuk,Gupta,icoz, Stancu, Sucu et al., Varma et al.}). Moreover,
important definitions, facts and features coupled with approximation theory
can be found in \cite{Altomare, Bernstein, DeVore-Lorentz, Korovkin,
Lorentz, Szasz}. For $\nu $,$\ x\in \lbrack 0,\infty )\ $and $g\in
C[0,\infty ),\ $in \cite{Sucu}, Dunkl analogue of Sz\'{a}sz operators is
given by 
\begin{equation}
S_{n}^{\ast }\left( g;x\right) =\frac{1}{e_{\nu }\left( nx\right) }%
\sum_{r=0}^{\infty }\frac{\left( nx\right) ^{r}}{\gamma _{\nu }\left(
r\right) }g\left( \frac{r+2\nu \theta _{r}}{n}\right) \ ;\ n\in 
\mathbb{N}
.  \label{1}
\end{equation}%
Here the $e_{\nu }\left( x\right) $ is defined as%
\begin{equation}
e_{\nu }\left( x\right) =\sum_{r=0}^{\infty }\frac{x^{r}}{\gamma _{\nu
}\left( r\right) },  \label{2}
\end{equation}%
where for $r\in 
\mathbb{N}
_{0}$ and $\nu >-\frac{1}{2}$ the coefficients $\gamma _{\nu }$ are given as
follows 
\begin{equation}
\gamma _{\nu }\left( 2r\right) =\frac{2^{2r}r!\Gamma \left( r+\nu
+1/2\right) }{\Gamma \left( \nu +1/2\right) }\text{ and }\gamma _{\nu
}\left( 2r+1\right) =\frac{2^{2r+1}r!\Gamma \left( r+\nu +3/2\right) }{%
\Gamma \left( \nu +1/2\right) }  \label{3}
\end{equation}%
in \cite{Rosenblum}. Also the coefficients $\gamma _{\nu }$ has the
recursion relation%
\begin{equation}
\frac{\gamma _{\nu }\left( r+1\right) }{\gamma _{\nu }\left( r\right) }%
=\left( 2\nu \theta _{r+1}+r+1\right) ,\text{ }r\in 
\mathbb{N}
_{0},  \label{4}
\end{equation}
where for $p\in 
\mathbb{N}
,\ \theta _{r}$ is given as 
\begin{equation}
\theta _{r}=\left\{ 
\begin{array}{cc}
0, & if\text{ }r=2p \\ 
1, & if\text{ }r=2p+1%
\end{array}%
\right. .  \label{5}
\end{equation}%
Also, the authors gave the other Dunkl generalizations of Sz\'{a}sz
operators in \cite{GB1, N}. Now, for $n\geq 1,\ $we define a Dunkl analogue
of Sz\'{a}sz beta type operators defined by%
\begin{equation}
T_{n}\left( g;x\right) =\frac{\left( n-1\right) }{e_{\nu }\left( nx\right) }%
\sum_{r=1}^{\infty }\tbinom{n+r-2}{r-1}\frac{\left( nx\right) ^{r}}{\gamma
_{\nu }\left( r\right) }\int_{0}^{\infty }s^{r-1}(1+s)^{-n-r+1}g\left(
s\right) ds+\frac{g(0)}{e_{\nu }\left( nx\right) },  \label{9}
\end{equation}%
where $e_{\nu }\left( x\right) $ and $\gamma _{\nu }$ are defined in $\left( %
\ref{2}\right) $ and $\left( \ref{3}\right) $, and $\nu ,x\geq 0$.
Furthermore $g\left( s\right) $ is defined on subset of all continuous
functions on $\left[ 0,\infty \right) $ for which the integral exists
finitely. Here well-known Beta function is denoted as $B(.,.)$ is given%
\begin{equation}
\int_{0}^{\infty }s^{r-1}(1+s)^{-n-r+1}ds=B(r,n-1).  \label{7}
\end{equation}

\begin{lemma}
\label{Lemma 1} Using (\ref{7}),\ we derive for $m\in 
\mathbb{N}
$%
\begin{equation}
\int_{0}^{\infty }s^{r-1}(1+s)^{-n-r+1}s^{m}ds=B(r+m,n-m-1).  \label{8}
\end{equation}
\end{lemma}

\begin{lemma}
\label{Lemma 2}For $T_{n}$ operators in (\ref{9}), the important properties
are hold.%
\begin{eqnarray}
T_{n}\left( 1;x\right) &=&1,\text{ }  \label{10} \\
&&  \notag \\
\left\vert T_{n}\left( s;x\right) -x\right\vert &\leq &\tfrac{2}{n-2}x+%
\tfrac{2\nu }{n-2}\text{ for}\ n>2,  \label{101}
\end{eqnarray}%
\begin{eqnarray}
\left\vert T_{n}\left( s^{2};x\right) -x^{2}\right\vert &\leq &\tfrac{5n-6}{%
n^{2}-5n+6}x^{2}+\left[ \tfrac{4\nu n}{n^{2}-5n+6}+\tfrac{2n}{n^{2}-5n+6}%
\right] x+\tfrac{4\nu ^{2}+6\nu }{n^{2}-5n+6}\text{ for}\ n>3,  \label{11} \\
&&  \notag \\
\left\vert T_{n}\left( s^{3};x\right) -x^{3}\right\vert &\leq &\tfrac{1}{%
(n-2)(n-3)(n-4)}\left\{ (9n^{2}-26n+24)x^{3}+6n^{2}(\nu +1)x^{2}\right. 
\notag \\
&&  \notag \\
&&\left. +(12\nu ^{2}n+18\nu n+6n)x+12\nu ^{2}+4\nu +8\nu ^{3}\right\} \text{
for}\ n>4,  \label{102} \\
&&  \notag \\
\left\vert T_{n}\left( s^{4};x\right) -x^{4}\right\vert &\leq &\tfrac{1}{%
n^{4}-14n^{3}+71n^{2}-154n+120}\left\{ (14n^{3}-71n^{2}+154n-120)x^{4}+(8\nu
n^{3}+12n^{3})x^{3}\right.  \notag \\
&&  \notag \\
&&+(62\nu n^{2}+24\nu ^{2}n^{2}+36n^{2})x^{2}+(32\nu ^{3}n+96\nu ^{2}n+56\nu
n+24n)x  \notag \\
&&  \notag \\
&&\left. +16\nu ^{4}+48\nu ^{3}+44\nu ^{2}+12\nu \right\} \text{ for}\ n>5.
\label{103}
\end{eqnarray}
\end{lemma}

\begin{lemma}
\label{Lemma 3}For $T_{n}$ operators, we have%
\begin{eqnarray}
\Psi _{1} &:&=T_{n}(s-x;x)\leq \left( \tfrac{2}{n-2}\right) x+\tfrac{2\nu }{%
n-2}\text{ for}\ n>2,  \notag \\
&&  \notag \\
\Psi _{2} &:&=T_{n}(\left( s-x\right) ^{2};x)\leq \left( \tfrac{9}{n-3}%
\right) x^{2}+\left( \tfrac{8\nu n-12\nu +2n}{n^{2}-5n+6}\right) x\text{ for}%
\ n>3,  \label{A} \\
&&  \notag \\
\Psi _{3} &:&=T_{n}((s-x)^{4};x)\leq \tfrac{88n^{2}-405n+450}{(n-2)(n-3)(n-5)%
}x^{4}+\tfrac{\left( 240n-480\nu +856n\nu -432n^{2}\nu +64n^{3}\nu
-228n^{2}+48n^{3}\right) \allowbreak }{n^{4}-14n^{3}+71n^{2}-154n+120}x^{3} 
\notag \\
&&  \notag \\
&&+\tfrac{\left( 96n^{2}\nu ^{2}+170n^{2}\nu +60n^{2}-456n\nu ^{2}-684n\nu
-120n+480\nu ^{2}+720\nu \right) }{n^{4}-14n^{3}+71n^{2}-154n+120}x^{2} 
\notag \\
&&  \notag \\
&&+\tfrac{\left( 24n-80\nu +72n\nu -240\nu ^{2}-160\nu ^{3}+144n\nu
^{2}+64n\nu ^{3}\right) }{n^{4}-14n^{3}+71n^{2}-154n+120}x  \notag \\
&&  \notag \\
&&+\tfrac{16\nu ^{4}+48\nu ^{3}+44\nu ^{2}+12\nu }{%
n^{4}-14n^{3}+71n^{2}-154n+120}\text{ for}\ n>6.  \label{B}
\end{eqnarray}
\end{lemma}

\begin{theorem}
\label{Theorem 4}For the operators in $\left( \ref{9}\right) $ and any $g\in
C[0,\infty )\cap E$, one obtain%
\begin{equation*}
T_{n}\left( g;x\right) \overset{\text{uniformly}}{\rightrightarrows }g\left(
x\right)
\end{equation*}%
on $A\subset $ $[0,\infty )$ which is each compact set\ as $n\rightarrow
\infty $. Here%
\begin{equation*}
E:=\{g:x\in \lbrack 0,\infty ),\frac{g\left( x\right) }{1+x^{2}}\text{ is
convergent as }x\rightarrow \infty \}\text{.}
\end{equation*}
\end{theorem}

Now, we evoke functions in the weighted spaces given on $[0,\infty )$ to
touch weighted approximation of our operators:%
\begin{eqnarray*}
B_{\rho }\left( 
\mathbb{R}
^{+}\right) &:&=\{g:\left\vert g\left( x\right) \right\vert \leq M_{g}\rho
\left( x\right) \}, \\
C_{\rho }\left( 
\mathbb{R}
^{+}\right) &:&=\{g:g\in C[0,\infty )\cap B_{\rho }\left( 
\mathbb{R}
^{+}\right) \}, \\
C_{\rho }^{\ast }\left( 
\mathbb{R}
^{+}\right) &:&=\{g:g\in C_{\rho }\left( 
\mathbb{R}
^{+}\right) \text{ and }\lim_{x\rightarrow \infty }\frac{g\left( x\right) }{%
\rho \left( x\right) }=cnst\}.
\end{eqnarray*}%
Here the weight function$\ $is called by $\rho \left( x\right) =1+x^{2}$ and 
$M_{g}$ is a constant based just on the function $g$. Also we keep in mind
the space $C_{\rho }\left( 
\mathbb{R}
^{+}\right) $ has a norm as $\left\Vert g\right\Vert _{\rho
}:=\sup\limits_{x\geq 0}\frac{\left\vert g\left( x\right) \right\vert }{\rho
\left( x\right) }$ (see \cite{Atakut}).

\begin{theorem}
\label{Theorem 5}For operators $T_{n}$ in $\left( \ref{9}\right) $ and each
function $g\in C_{\rho }^{\ast }\left( 
\mathbb{R}
^{+}\right) ,$ one has%
\begin{equation*}
\lim_{n\rightarrow \infty }\left\Vert T_{n}\left( g;x\right) -g\left(
x\right) \right\Vert _{\rho }=0.
\end{equation*}
\end{theorem}

\section{Convergence of operators in (\protect\ref{9})}

Firstly, we remind the Lipschitz class of order $\alpha $ for function $g$.
If $g\in Lip_{M}\left( \alpha \right) $,\ then $g$ satisfies the inequality%
\begin{equation*}
\left\vert g\left( s\right) -g\left( t\right) \right\vert \leq M\left\vert
s-t\right\vert ^{\alpha }
\end{equation*}%
where $s,t\in \lbrack 0,\infty ),\ 0<\alpha \leq 1$ and $M>0.$

\begin{theorem}
\label{Theorem 6}If $h\in Lip_{M}\left( \alpha \right) $, the following
inequality%
\begin{equation*}
\left\vert T_{n}\left( h;x\right) -h\left( x\right) \right\vert \leq M\left(
\tau _{n}\left( x\right) \right) ^{\alpha /2},
\end{equation*}%
is hold where $\tau _{n}\left( x\right) =\Psi _{2}.$
\end{theorem}

Now, we deal with the space symbolized by $\overset{\sim }{C}[0,\infty )$
has uniformly continuous functions on $[0,\infty )$ and modulus of
continuity $g\in \overset{\sim }{C}[0,\infty )$ is denoted as%
\begin{equation}
\omega \left( g;\delta \right) :=\sup\limits_{\substack{ s,t\in \lbrack
0,\infty )  \\ \left\vert s-t\right\vert \leq \delta }}\left\vert g\left(
s\right) -g\left( t\right) \right\vert .  \label{12}
\end{equation}

\begin{theorem}
\label{Theorem 7}The operators in (\ref{9}) satisfy the inequality%
\begin{equation*}
\left\vert T_{n}\left( g;x\right) -g\left( x\right) \right\vert \leq \left(
1+\sqrt{9x^{2}+\left( \tfrac{8\nu n-12\nu +2n}{n-2}\right) x}\right) \omega
\left( g;\frac{1}{\sqrt{n}}\right) ,
\end{equation*}%
where $g\in \widetilde{C}[0,\infty )\cap E$ ,$\ \omega $ is modulus of
continuity.
\end{theorem}

Now, we note that the space $C_{B}[0,\infty )$ is all continuous and bounded
functions on $[0,\infty )$. Also%
\begin{equation}
C_{B}^{2}[0,\infty )=\{g\in C_{B}[0,\infty ):g^{\prime },g^{\prime \prime
}\in C_{B}[0,\infty )\}  \label{16}
\end{equation}%
and the norm on $C_{B}^{2}[0,\infty )$ is defined as%
\begin{equation*}
\left\Vert g\right\Vert _{C_{B}^{2}[0,\infty )}=\left\Vert g\right\Vert
_{C_{B}[0,\infty )}+\left\Vert g^{\prime }\right\Vert _{C_{B}[0,\infty
)}+\left\Vert g^{\prime \prime }\right\Vert _{C_{B}[0,\infty )}
\end{equation*}%
for $\forall g\in C_{B}^{2}[0,\infty ).$

\begin{lemma}
\label{Lemma 8}For $h\in C_{B}^{2}[0,\infty )$, one has the inequality%
\begin{equation}
\left\vert T_{n}\left( h;x\right) -h\left( x\right) \right\vert \leq \chi
_{n}\left( x\right) \left\Vert h\right\Vert _{C_{B}^{2}[0,\infty )},
\label{17}
\end{equation}%
where%
\begin{equation}
\chi _{n}\left( x\right) =\Psi _{1}+\Psi _{2}.  \label{18}
\end{equation}
\end{lemma}

Note that the second order of modulus continuity of $g$ on $C_{B}[0,\infty )$
is as%
\begin{equation*}
\omega _{2}\left( g;\delta \right) :=\sup_{0<s\leq \delta }\left\Vert
g\left( .+2s\right) -2g\left( .+s\right) +g\left( .\right) \right\Vert
_{C_{B}[0,\infty )}.
\end{equation*}%
Thus, we can obtain the following important theorem.

\begin{theorem}
\label{Theorem 9}The operators in (\ref{9}) satisfy the following inequality%
\begin{equation}
\left\vert T_{n}\left( g;x\right) -g\left( x\right) \right\vert \leq
2M\left\{ \min \left( 1,\frac{\chi _{n}\left( x\right) }{2}\right)
\left\Vert g\right\Vert _{C_{B}[0,\infty )}+\omega _{2}\left( g;\sqrt{\frac{%
\chi _{n}\left( x\right) }{2}}\right) \right\}  \label{20}
\end{equation}%
where $\forall g\in C_{B}[0,\infty ),\ x\in \lbrack 0,\infty )$, $M$ is a
positive constant which is not based on $n$ and $\chi _{n}\left( x\right) $
is in $\left( \ref{18}\right) .$
\end{theorem}

Now, we focus the order of the functions $g\in C_{\rho }^{\ast }(%
\mathbb{R}
^{+}).$ Atakut and Ispir \cite{Atakut}, Ispir \cite{ispir} defined the
weighted of continuity denoted by%
\begin{equation*}
\Omega (g;\delta )=\sup_{x\in \lbrack 0,\infty ),\ \left\vert h\right\vert
\leq \delta }\frac{\left\vert g(x+h)-g(x)\right\vert }{(1+h^{2})(1+x^{2})}
\end{equation*}%
for $g\in C_{\rho }^{\ast }(%
\mathbb{R}
^{+}).$ This modulus satisfying $\lim\limits_{\delta \rightarrow 0}\Omega
(g;\delta )=0$ and%
\begin{equation}
\left\vert g(s)-g(x)\right\vert \leq 2\left( 1+\frac{\left\vert
s-x\right\vert }{\delta }\right) (1+\delta ^{2})(1+x^{2})(1+(s-x)^{2})\Omega
(g;\delta ),  \label{21}
\end{equation}%
where $s,x\in \lbrack 0,\infty ).$

\begin{theorem}
\label{Theorem 10}The operators in (\ref{9}) satisfy the following inequality%
\begin{equation*}
\sup_{x\in \lbrack 0,\infty )}\frac{\left\vert T_{n}(g;x)-g(x)\right\vert }{%
(1+x^{2})^{3}}\leq M_{\nu }\left( 1+\frac{1}{n}\right) \Omega \left( g;\frac{%
1}{\sqrt{n}}\right) ,
\end{equation*}%
where $g\in C_{\rho }^{\ast }(%
\mathbb{R}
^{+})$ and $M_{\nu }$ is a constant which is not based on $x.$%
\begin{equation*}
\end{equation*}
\end{theorem}

\section{The proofs of the results}

\textbf{Proof of Lemma \ref{Lemma 2}.\ }

For $g\left( s\right) =1,$ using (\ref{8}) and $e_{\nu }\left( x\right) ,$
we have%
\begin{eqnarray*}
T_{n}\left( 1;x\right) &=&\frac{\left( n-1\right) }{e_{\nu }\left( nx\right) 
}\sum_{r=1}^{\infty }\frac{\left( nx\right) ^{r}}{\gamma _{\nu }\left(
r\right) }\tbinom{n+r-2}{r-1}\int_{0}^{\infty }\frac{s^{r-1}}{(1+s)^{n+r-1}}%
ds+\frac{1}{e_{\nu }\left( nx\right) } \\
&=&\frac{1}{e_{\nu }\left( nx\right) }\sum_{r=1}^{\infty }\frac{\left(
nx\right) ^{r}}{\gamma _{\nu }\left( r\right) }+\frac{1}{e_{\nu }\left(
nx\right) }=1.
\end{eqnarray*}%
For $n>2$ and $g\left( s\right) =s,$ using (\ref{8}), (\ref{4}) and $e_{\nu
}\left( x\right) ,\ $respectively,$\ $one derive%
\begin{eqnarray*}
T_{n}\left( s;x\right) &=&\frac{n-1}{e_{\nu }\left( nx\right) }%
\sum_{r=1}^{\infty }\frac{\left( nx\right) ^{r}}{\gamma _{\nu }\left(
r\right) }\int_{0}^{\infty }\frac{s^{r}}{(1+s)^{n+r-1}}\tbinom{n+r-2}{r-1}ds+%
\frac{g(0)}{e_{\nu }\left( nx\right) } \\
&=&\frac{1}{e_{\nu }\left( nx\right) }\sum_{r=1}^{\infty }\frac{\left(
nx\right) ^{r}}{\gamma _{\nu }\left( r\right) }\frac{(r+2\nu \theta
_{r}-2\nu \theta _{r})}{n-2} \\
&=&\frac{1}{e_{\nu }\left( nx\right) }\frac{1}{n-2}\sum_{r=1}^{\infty }\frac{%
\left( nx\right) ^{r}}{\gamma _{\nu }\left( r-1\right) }-\frac{1}{e_{\nu
}\left( nx\right) }\frac{2\nu }{n-2}\sum_{r=1}^{\infty }\theta _{r}\frac{%
\left( nx\right) ^{r}}{\gamma _{\nu }\left( r\right) } \\
&=&\frac{1}{n-2}\left( nx-\frac{2\nu }{e_{\nu }\left( nx\right) }%
\sum_{r=1}^{\infty }\theta _{r}\frac{\left( nx\right) ^{r}}{\gamma _{\nu
}\left( r\right) }\right) .
\end{eqnarray*}%
Thus, we derive%
\begin{eqnarray*}
\left\vert T_{n}\left( s;x\right) -x\right\vert &\leq &\frac{1}{n-2}\left(
2x+\frac{2\nu }{e_{\nu }\left( nx\right) }\sum_{r=0}^{\infty }\frac{\left(
nx\right) ^{r}}{\gamma _{\nu }\left( r\right) }\right) \\
&\leq &\left( \frac{2}{n-2}\right) x+\frac{2\nu }{n-2}.
\end{eqnarray*}%
For $n>3$ and $g\left( s\right) =s^{2},$ using (\ref{8}),$\ $we get%
\begin{eqnarray*}
T_{n}\left( s^{2};x\right) &=&\frac{n-1}{e_{\nu }\left( nx\right) }%
\sum_{r=1}^{\infty }\frac{\left( nx\right) ^{r}}{\gamma _{\nu }\left(
r\right) }\int_{0}^{\infty }\frac{s^{r+1}}{(1+s)^{n+r-1}}\tbinom{n+r-2}{r-1}%
ds+\frac{g(0)}{e_{\nu }\left( nx\right) } \\
&=&\frac{1}{n^{2}-5n+6}\frac{1}{e_{\nu }\left( nx\right) }\sum_{r=1}^{\infty
}r(r+1)\frac{\left( nx\right) ^{r}}{\gamma _{\nu }\left( r\right) }.
\end{eqnarray*}%
Using%
\begin{equation*}
r(r+1)=-(r-1)2\nu \theta _{r}+(r+2\nu \theta _{r})(r-1)+2r
\end{equation*}%
and (\ref{4}), we have%
\begin{equation*}
\begin{array}{l}
T_{n}\left( s^{2};x\right) =\frac{1}{(n^{2}-5n+6)}\frac{1}{e_{\nu }\left(
nx\right) }\left\{ \sum_{r=1}^{\infty }\frac{\left( nx\right) ^{r}}{\gamma
_{\nu }\left( r\right) }(r+2\nu \theta _{r})(r-1)\right. \\ 
\\ 
-\left. 2\nu \sum_{r=1}^{\infty }\theta _{r}(r-1)\frac{\left( nx\right) ^{r}%
}{\gamma _{\nu }\left( r\right) }+\sum_{r=1}^{\infty }\frac{\left( nx\right)
^{r}}{\gamma _{\nu }\left( r\right) }2r\right\} \\ 
\\ 
=\frac{1}{(n-2)(n-3)}\frac{1}{e_{\nu }\left( nx\right) }\left\{
\sum_{r=0}^{\infty }(r+2\nu \theta _{r}-2\nu \theta _{r})\frac{\left(
nx\right) ^{r+1}}{\gamma _{\nu }\left( r\right) }\right. \\ 
\\ 
\left. -2\nu \sum_{r=1}^{\infty }\theta _{r}(r-1)\frac{\left( nx\right) ^{r}%
}{\gamma _{\nu }\left( r\right) }+\sum_{r=1}^{\infty }\frac{\left( nx\right)
^{r}}{\gamma _{\nu }\left( r\right) }2r\right\} \\ 
\\ 
=\frac{1}{n^{2}-5n+6}\frac{1}{e_{\nu }\left( nx\right) }\left\{
\sum_{r=1}^{\infty }\frac{\left( nx\right) ^{r+1}}{\gamma _{\nu }\left(
r-1\right) }-2\nu \sum_{r=1}^{\infty }\theta _{r}\frac{\left( nx\right)
^{r+1}}{\gamma _{\nu }\left( r\right) }\right. \\ 
\\ 
\left. -2\nu \sum_{r=1}^{\infty }\theta _{r}(r-1)\frac{\left( nx\right) ^{r}%
}{\gamma _{\nu }\left( r\right) }+\sum_{r=1}^{\infty }\frac{\left( nx\right)
^{r}}{\gamma _{\nu }\left( r\right) }2r\right\} .%
\end{array}%
\end{equation*}%
Therefore, since\ $\theta _{r}\leq 1,$ we obtain%
\begin{equation*}
\begin{array}{l}
\left\vert T_{n}\left( s^{2};x\right) -x^{2}\right\vert \leq \frac{1}{\left(
n^{2}-5n+6\right) e_{\nu }\left( nx\right) }\left\{ \left\vert
\sum_{r=1}^{\infty }\frac{\left( nx\right) ^{r+1}}{\gamma _{\nu }\left(
r-1\right) }-(n-2)(n-3)e_{\nu }\left( nx\right) x^{2}\right\vert \right. \\ 
\\ 
+\left. 2\nu \sum_{r=1}^{\infty }\frac{\left( nx\right) ^{r+1}}{\gamma _{\nu
}\left( r\right) }+2\nu \sum_{r=1}^{\infty }(r-1)\frac{\left( nx\right) ^{r}%
}{\gamma _{\nu }\left( r\right) }+\sum_{r=1}^{\infty }\frac{\left( nx\right)
^{r}}{\gamma _{\nu }\left( r\right) }2r\right\} \\ 
\\ 
\leq \frac{1}{(n-2)(n-3)}\frac{1}{e_{\nu }\left( nx\right) }\left\{
\left\vert n^{2}x^{2}e_{\nu }\left( nx\right) -(n-2)(n-3)e_{\nu }\left(
nx\right) x^{2}\right\vert +2\nu nxe_{\nu }\left( nx\right) \right. \\ 
\\ 
\left. +2\nu \sum_{r=1}^{\infty }\frac{\left( nx\right) ^{r}}{\gamma _{\nu
}\left( r\right) }(r+2\nu \theta _{r}-2\nu \theta
_{r}-1)+2\sum_{r=1}^{\infty }\frac{\left( nx\right) ^{r}}{\gamma _{\nu
}\left( r\right) }\left( r+2\nu \theta _{r}-2\nu \theta _{r}\right) \right\}
\\ 
\\ 
\leq \frac{1}{n^{2}-5n+6}\frac{1}{e_{\nu }\left( nx\right) }\left\{
\left\vert n^{2}x^{2}e_{\nu }\left( nx\right) -(n-2)(n-3)e_{\nu }\left(
nx\right) x^{2}\right\vert \right. \\ 
\\ 
+\left. 2\nu nxe_{\nu }\left( nx\right) +2\nu (nx+(2\nu +1))e_{\nu }\left(
nx\right) +2(nx+2\nu )e_{\nu }\left( nx\right) \right\} \\ 
\\ 
=\frac{5n-6}{n^{2}-5n+6}x^{2}+\left[ \frac{4\nu n}{n^{2}-5n+6}+\frac{2n}{%
n^{2}-5n+6}\right] x+\frac{4\nu ^{2}+6\nu }{n^{2}-5n+6}.%
\end{array}%
\end{equation*}%
Similarly, (\ref{102}) and (\ref{103}) can be proved.

\textbf{Proof of Theorem \ref{Theorem 4}. }As $n\rightarrow \infty ,\ $under
favour of Korovkin Theorem in \cite{Korovkin}, one has $T_{n}\left(
g;x\right) \overset{\text{uniformly}}{\rightrightarrows }g\left( x\right) $
on $A\subset \lbrack 0,\infty )$ which is each compact set because $%
\lim_{n\rightarrow \infty }T_{n}(e_{i};x)=x^{i},$ for$\ e_{i}=s^{i},\
i=0,1,2,$ which is uniformly on $A\subset \lbrack 0,\infty )$ with the help
of using Lemma \textbf{\ref{Lemma 2}}.

\textbf{Proof of Theorem \ref{Theorem 5}. }From $(\ref{10})$, we can write $%
\lim_{n\rightarrow \infty }\left\Vert T_{n}\left( 1;x\right) -1\right\Vert
_{\rho }=0.$

For$\ n>2,\ $by (\ref{101}) and the following calculation%
\begin{eqnarray*}
\sup_{x\in \lbrack 0,\infty )}\frac{\left\vert T_{n}\left( s;x\right)
-x\right\vert }{1+x^{2}} &\leq &\frac{1}{n-2}\left( 2\sup_{x\in \lbrack
0,\infty )}\frac{x}{1+x^{2}}+2\nu \sup_{x\in \lbrack 0,\infty )}\frac{1}{%
1+x^{2}}\right) \\
&\leq &\left( \frac{2}{n-2}\right) +\frac{2\nu }{n-2},
\end{eqnarray*}%
we get%
\begin{equation*}
\lim_{n\rightarrow \infty }\left\Vert T_{n}\left( s;x\right) -x\right\Vert
_{\rho }=0.
\end{equation*}%
Finally, for$\ n>3,$ by $(\ref{11})$ and the following calculation%
\begin{eqnarray*}
\sup_{x\in \lbrack 0,\infty )}\frac{\left\vert T_{n}\left( s^{2};x\right)
-x^{2}\right\vert }{1+x^{2}} &\leq &\frac{1}{(n^{2}-5n+6)}\left\{ \left(
5n-6\right) \sup_{x\in \lbrack 0,\infty )}\frac{x^{2}}{1+x^{2}}\right. \\
&&+\left. \left[ 4\nu n+2n\right] \sup_{x\in \lbrack 0,\infty )}\frac{x}{%
1+x^{2}}\right. \\
&&+\left. \left( 4\nu ^{2}+6\nu \right) \sup_{x\in \lbrack 0,\infty )}\frac{1%
}{1+x^{2}}\right\} \\
&=&\frac{1}{(n^{2}-5n+6)}\left\{ \left( 5n-6\right) +\left( 2\nu n+n\right)
+\left( 4\nu ^{2}+6\nu \right) \right\} ,
\end{eqnarray*}%
we have%
\begin{equation*}
\lim_{n\rightarrow \infty }\left\Vert T_{n}\left( s^{2};x\right)
-x^{2}\right\Vert _{\rho }=0.
\end{equation*}%
Thus, we get $\lim_{n\rightarrow \infty }\left\Vert T_{n}\left( g;x\right)
-g\left( x\right) \right\Vert _{\rho }=0$ for each $g\in C_{\rho }^{\ast
}\left( 
\mathbb{R}
^{+}\right) $ via weighted Korovkin-type theorem given by Gadzhiev \cite%
{Gadzhiev}.

\textbf{Proof of Theorem \ref{Theorem 6}.} Using $h\in Lip_{M}\left( \alpha
\right) $ and linearity, one has%
\begin{eqnarray*}
\left\vert T_{n}\left( h;x\right) -h\left( x\right) \right\vert &\leq
&T_{n}\left( \left\vert h\left( s\right) -h\left( x\right) \right\vert
;x\right) \\
&\leq &MT_{n}\left( \left\vert s-x\right\vert ^{\alpha };x\right) .
\end{eqnarray*}%
From Lemma \textbf{\ref{Lemma 2}} and H\"{o}lder's famous inequality, we
derive%
\begin{equation*}
\left\vert T_{n}\left( h;x\right) -h\left( x\right) \right\vert \leq M\left[
\Psi _{2}\right] ^{\frac{\alpha }{2}}.
\end{equation*}%
Then choosing $\tau _{n}\left( x\right) =\Psi _{2}$, thus one has the
required inequality.

\textbf{Proof of Theorem \ref{Theorem 7}.} By the property of modulus of
continuity and (\ref{A}), one get%
\begin{eqnarray*}
\left\vert T_{n}\left( g;x\right) -g\left( x\right) \right\vert &\leq
&T_{n}\left( \left\vert g\left( s\right) -g\left( x\right) \right\vert
;x\right) \\
&\leq &\left( 1+\frac{1}{\delta }T_{n}\left( \left\vert s-x\right\vert
;x\right) \right) \omega \left( g;\delta \right) \\
&\leq &\left( 1+\frac{1}{\delta }\sqrt{\Psi _{2}}\right) \omega \left(
g;\delta \right) .
\end{eqnarray*}%
Then using Cauchy-Schwarz's famous inequality, one has%
\begin{equation}
\left\vert T_{n}\left( g;x\right) -g\left( x\right) \right\vert \leq \left(
1+\frac{1}{\delta }\sqrt{\tfrac{9}{n-3}x^{2}+\left( \tfrac{8\nu n-12\nu +2n}{%
n^{2}-5n+6}\right) x}\right) \omega \left( g;\delta \right) .  \label{15}
\end{equation}%
Choosing $\delta =\frac{1}{\sqrt{n}},$ the proof is done.

\textbf{Proof of Lemma \ref{Lemma 8}. }Using the Taylor's series of the
function $h$, we can write%
\begin{equation*}
h\left( s\right) =h\left( x\right) +\left( s-x\right) h^{\prime }\left(
x\right) +\frac{\left( s-x\right) ^{2}}{2!}h^{\prime \prime }\left( \varrho
\right) ,\text{ }\varrho \in \left( x,s\right) .
\end{equation*}%
From the linear operator, we give%
\begin{equation*}
T_{n}\left( h;x\right) -h\left( x\right) =h^{\prime }\left( x\right) \Psi
_{1}+\frac{h^{\prime \prime }\left( \varrho \right) }{2}\Psi _{2}.
\end{equation*}%
Then for $n>3,$ using Lemma \textbf{\ref{Lemma 3}}, one obtain%
\begin{eqnarray*}
\left\vert T_{n}\left( h;x\right) -h\left( x\right) \right\vert &\leq &\frac{%
2}{n-2}\left( x+\nu \right) \left\Vert h^{\prime }\right\Vert
_{C_{B}[0,\infty )} \\
&& \\
&&+\left[ \left( \tfrac{9}{n-3}\right) x^{2}+\left( \tfrac{8\nu n-12\nu +2n}{%
n^{2}-5n+6}\right) x\right] \left\Vert h^{\prime \prime }\right\Vert
_{C_{B}[0,\infty )} \\
&& \\
&\leq &[\Psi _{1}+\Psi _{2}]\left\Vert h\right\Vert _{C_{B}^{2}[0,\infty )}.
\end{eqnarray*}

Choosing $\chi _{n}\left( x\right) =[\Psi _{1}+\Psi _{2}]$ which finishes
the proof.

\textbf{Proof of Theorem \ref{Theorem 9}. }For any $f\in C_{B}^{2}[0,\infty
) $, from the triangle inequality and Lemma \textbf{\ref{Lemma 8}}, one has%
\begin{eqnarray*}
\Lambda &=&\left\vert T_{n}\left( g;x\right) -g\left( x\right) \right\vert
\leq \left\vert T_{n}\left( g-f;x\right) \right\vert +\left\vert T_{n}\left(
f;x\right) -f\left( x\right) \right\vert +\left\vert g\left( x\right)
-f\left( x\right) \right\vert \\
&& \\
&\leq &2\left\Vert g-f\right\Vert _{C_{B}[0,\infty )}+\chi _{n}\left(
x\right) \left\Vert f\right\Vert _{C_{B}^{2}[0,\infty )} \\
&& \\
&=&2\left\{ \left\Vert g-f\right\Vert _{C_{B}[0,\infty )}+\frac{\chi _{n}}{2}%
\left( x\right) \left\Vert f\right\Vert _{C_{B}^{2}[0,\infty )}\right\} .
\end{eqnarray*}%
With the help of Peetre's $K$ functional in \cite{DeVore-Lorentz}, one has 
\begin{equation*}
\Lambda \leq 2K\left( g;\frac{\chi _{n}\left( x\right) }{2}\right) .
\end{equation*}%
Thus we can write%
\begin{equation*}
\Lambda \leq 2M\left\{ \min \left( 1,\frac{\chi _{n}\left( x\right) }{2}%
\right) \left\Vert g\right\Vert _{C_{B}[0,\infty )}+\omega _{2}\left( g;%
\sqrt{\frac{\chi _{n}\left( x\right) }{2}}\right) \right\}
\end{equation*}%
because of the well-known connection between $K_{2}$ and $\omega _{2}$ in 
\cite{DeVore-Lorentz}. We note that the connection is as%
\begin{equation*}
K_{2}(g;\delta )\leq C\left\{ \min (1,\delta )\left\Vert g\right\Vert
_{C_{B}[0,\infty )}+\omega _{2}\left( g;\sqrt{\delta }\right) \right\} .
\end{equation*}%
Here $C$ is an positive constant \cite{Ciupa}.

\textbf{Proof of Theorem \ref{Theorem 10}.}

Under favour of (\ref{21}) and the property of linearity of operator, one has%
\begin{eqnarray*}
\Lambda &=&\left\vert T_{n}(g;x)-g(x)\right\vert \leq T_{n}(\left\vert
g(s)-g(x)\right\vert ;x) \\
&& \\
&\leq &2(1+\delta ^{2})(1+x^{2})\Omega (g;\delta )T_{n}\left( \left( 1+\frac{%
\left\vert s-x\right\vert }{\delta }\right) (1+(s-x)^{2});x\right) \\
&& \\
&=&2(1+\delta ^{2})(1+x^{2})\Omega (g;\delta )\left\{ T_{n}\left( 1;x\right)
+\frac{1}{\delta }T_{n}\left( \left\vert s-x\right\vert ;x\right) \right. \\
&&\left. +T_{n}\left( \left( s-x\right) ^{2};x\right) +\frac{1}{\delta }%
T_{n}\left( \left\vert s-x\right\vert ^{3};x\right) \right\} .
\end{eqnarray*}

Now, if we apply Cauchy-Scwarz's inequality for $T_{n}\left( \left\vert
s-x\right\vert ;x\right) $ and $T_{n}\left( \left\vert s-x\right\vert
^{3};x\right) ,$ then we derive%
\begin{eqnarray*}
T_{n}\left( \left\vert s-x\right\vert ;x\right) &\leq &\sqrt{T_{n}\left(
\left( s-x\right) ^{2};x\right) }, \\
T_{n}\left( \left\vert s-x\right\vert ^{3};x\right) &\leq &\sqrt{T_{n}\left(
\left( s-x\right) ^{2};x\right) }\sqrt{T_{n}\left( \left( s-x\right)
^{4};x\right) }.
\end{eqnarray*}%
Thus, we get 
\begin{equation*}
\Lambda \leq 2(1+\delta ^{2})(1+x^{2})\Omega (g;\delta )\left( 1+\frac{1}{%
\delta }\sqrt{\Psi _{2}}+\Psi _{2}+\frac{1}{\delta }\sqrt{\Psi _{2}\Psi _{3}}%
\right) .
\end{equation*}%
With the help of $\left( \ref{A}\right) $ and $\left( \ref{B}\right) $, one
has%
\begin{eqnarray*}
\Lambda &\leq &2(1+\delta ^{2})(1+x^{2})\Omega (g;\delta )\left\{ 1+\tfrac{9%
}{n-3}x^{2}+\left( \tfrac{8\nu n-12\nu +2n}{n^{2}-5n+6}\right) x\right. \\
&& \\
&&+\frac{1}{\delta }\sqrt{\tfrac{9}{n-3}x^{2}+\left( \tfrac{8\nu n-12\nu +2n%
}{n^{2}-5n+6}\right) x} \\
&& \\
&&\left. +\frac{1}{\delta }\sqrt{\left( \tfrac{9}{n-3}x^{2}+\left( \tfrac{%
8\nu n-12\nu +2n}{n^{2}-5n+6}\right) x\right) \Psi _{3}}\right\} .
\end{eqnarray*}%
Choosing $\delta =\frac{1}{\sqrt{n}},$ then the proof is completed.

\textbf{Acknowledgement. }The authors are grateful to the referees for their
valuable comments and suggestions which improved the quality and the clarity
of the paper.


\begin{thebibliography}{99}
\bibitem{Altomare} Altomare, F., Campiti, M. Korovkin-Type Approximation
Theory and its Applications, de Gruyter Studies in Mathematics, vol. 17,
Walter de Gruyter, Berlin, Germany, \textbf{1994}.

\bibitem{Atakut} Atakut, \c{C}., \.{I}spir, N. Approximation by modified Sz%
\'{a}sz--Mirakjan operators on weighted spaces, Proc. Indian Acad. Sci.
Math. 112 \textbf{(2002)}, 571--578

\bibitem{Ata-Buyuk} Atakut, \c{C}., B\"{u}y\"{u}kyazici, \.{I}. Stancu type
generalization of the Favard Sz\'{a}sz operators, Appl. Math. Lett., 23 (12) 
\textbf{(2010)}, 1479-1482.

\bibitem{Bernstein} Bernstein, S.N. D\'{e}monstration du th\'{e}or\'{e}me de
Weierstrass fond\'{e}e sur le calcul des probabilit\'{e}s, Commun. Soc.
Math. Kharkow 2 (13)\textbf{\ (1912)}, 1--2.

\bibitem{Ciupa} Ciupa, A. A class of integral Favard--Sz\'{a}sz type
operators. Stud. Univ. Babes-Bolyai Math. 40 (1) \textbf{(1995)}, 39--47.

\bibitem{DeVore-Lorentz} DeVore, R.A., Lorentz, G.G. Construtive
Approximation, Springer, Berlin, \textbf{1993}.

\bibitem{Gadzhiev} Gadzhiev, A.D. The convergence problem for a sequence of
positive linear operators on unbounded sets and theorems analogues to that
of P.P. Korovkin, Sov. Math. Dokl. 15 (5) \textbf{(1974)}, 1453-1436.

\bibitem{Gupta} Gupta, V., Vasishtha, V., Gupta, M.K. Rate of convergence of
the Sz\'{a}sz--Kantorovich--Bezier operators for bounded variation
functions, Publ. Inst. Math. (Beograd) (N.S.), 72 \textbf{(2006)}, 137--143.

\bibitem{GB1} \.{I}\c{c}\"{o}z, G., \c{C}ekim, B. Dunkl generalization of sz%
\'{a}sz operators via q-calculus. Journal of Inequalities and Applications
2015; 2015(Article ID 284):1--11.30.

\bibitem{icoz} \.{I}\c{c}\"{o}z, G., \c{C}ekim, B. Stancu-type
generalization of Dunkl analogue of Sz\'{a}sz--Kantorovich operators, Math.
Meth. Appl. Sci, 39 \textbf{(2016)}, 1803--1810.

\bibitem{ispir} \.{I}spir, N., On Modified Baskakov Operators on Weighted
Spaces, Turk. J. Math. 25 \textbf{(2001)}, 355 -- 365.

\bibitem{Korovkin} Korovkin, P. P. On convergence of linear positive
operators in the space of continuous functions (Russian), Doklady Akad.
Nauk. SSSR (NS) 90 \textbf{(1953)}, 961--964.

\bibitem{Lorentz} Lorentz,G.G. Bernstein polynomials, In Mathematical
Expositions, vol. 8, University of Toronto Press, Toronto, \textbf{1953}.

\bibitem{N} Wafi, A., Rao, N. Sz\'{a}sz-Durrmeyer Operators Based on Dunkl
Analogue. Complex Analysis and Operator Theory, DOI:
10.1007/s11785-017-0647-7.

\bibitem{Rosenblum} Rosenblum, M. Generalized Hermite polynomials and the
Bose-like oscillator calculus, Oper. Theory: Adv. Appl. 73 \textbf{(1994)},
369-396.

\bibitem{Stancu} Stancu, D.D. Approximation of function by a new class of
polynomial operators, Rev. Rourn. Math. Pures et Appl., 13 (8) \textbf{(1968)%
}, 1173-1194.

\bibitem{Sucu} Sucu, S. Dunkl analogue of Sz\'{a}sz operators, Appl. Math.
Comput. 244 \textbf{(2014)}, 42-48.

\bibitem{Sucu et al.} Sucu, S., \.{I}\c{c}\"{o}z, G., Varma, S. On some
extensions of Szasz operators including Boas-Buck type polynomials, Abstr.
Appl. Anal., Vol. \textbf{2012}, Article ID 680340, 15 pages.

\bibitem{Szasz} Sz\'{a}sz, O. Generalization of S. Bernstein's polynomials
to the infinite interval, J. Res. Nat. Bur. Stand. 45 \textbf{(1950)},
239--245.

\bibitem{Varma et al.} Varma, S., Sucu, S., \.{I}\c{c}\"{o}z, G.
Generalization of Szasz operators involving Brenke type polynomials, Comput.
Math. Appl., 64 (2)\textbf{\ (2012)}, 121-127.
\end{thebibliography}
\end{document}